\documentclass[12pt, oneside,reqno]{amsart}   	
\usepackage{geometry}
\usepackage{float}
\allowdisplaybreaks
\usepackage{graphicx}				
\usepackage{amssymb}
\usepackage{amsmath}
\usepackage[all]{xy}
\usepackage{mathtools}
\usepackage{tikz-cd}
\usepackage{enumitem}
\usepackage{tikz-cd}
\usepackage{appendix}

\usepackage[%
bookmarks=true,
colorlinks,
linkcolor=blue,
urlcolor=blue,
citecolor=blue,
plainpages=false,
pdfpagelabels,
final,
breaklinks=true\between
]{hyperref}
\usepackage{hyperref}
\usepackage[numbers,sort&compress]{natbib}

\newtheorem{thm}{Theorem}[section]
\newtheorem{lem}[thm]{Lemma}

\numberwithin{equation}{section}

\def\Pb{\ifmmode{\Bbb P}\else{$\Bbb P$}\fi}
\def\Z{\ifmmode{\Bbb Z}\else{$\Bbb Z$}\fi}
\def\C{\ifmmode{\Bbb C}\else{$\Bbb C$}\fi}
\def\R{\ifmmode{\Bbb R}\else{$\Bbb R$}\fi}
\def\S{\ifmmode{S^2}\else{$S^2$}\fi}

\def\S{\cal S}

\newenvironment{pf}{\paragraph{Proof:}}{\hfill$\square$ \newline}

\begin{document}

\title[medium entropy]{On self shrinkers of medium entropy in $\R^4$}
\author{Alexander Mramor}
\address{Department of Mathematics, Johns Hopkins University, Baltimore, MD, 21231}
\email{amramor1@jhu.edu}

\begin{abstract} In this article we study smooth asymptotically conical self shrinkers in $\R^4$ with Colding-Minicozzi entropy bounded above by $\Lambda_{1}$. 
\end{abstract}
\maketitle
\section{Introduction}

Self shrinkers are basic singularity models for the mean curvature flow and in the noncompact case nongeneric ones (generic ones being generalized round cylinders $S^k(\sqrt{2k}) \times \R^{n-k}$) are expected to often be asymptotically conical. The purpose of this paper is to understand the topology of smooth such self shrinkers $M^3 \subset \R^4$ satisfying their Colding-Minicozzi entropy $\lambda(M)$, discussed in section 2, bounded above by $\Lambda_1$, the entropy of the round circle.
$\medskip$

The main result of the paper is the following, in part inspired by arguments in \cite{BW, IW, MW, Mra1, W3, HW}. The basic idea is by considering renormalized mean curvature flows out of (appropriate perturbations of) asymptotically conical self shrinkers below we may use the entropy assumption to constrain which types of singularities may occur. This has strong implications for how topology may change under the flow. This is useful because topology can in a sense be used to ``trap'' the flow but on the other hand the flow must clear out; these two principles can then be played together to constrain the topology of the self shrinker in question.  
\begin{thm}\label{thm1} Suppose $M^3 \subset \R^{4}$ is a smooth 2-sided asymptotically conical self shrinker with entropy less than $\Lambda_{1}$ and $k$ ends. Then it is diffeomorphic to $S^3$ with $k$ 3-balls removed and replaced with $k$ copies of $S^{2} \times \R_+$ attached along their respective boundaries. If $k = 1$ then $M \simeq \R^3$ and in particular this is the case when $\lambda(M) \leq \Lambda_2$. 
\end{thm} 
 This extends to the noncompact case joint work of the author and S. Wang \cite{MW} on compact self shrinkers $M^n \subset \R^{n+1}$ when $n = 3$, where they showed for each $n \geq 3$ closed self shrinkers $M^n$ with $\Lambda(M) < \Lambda_{n-2}$ are diffeomorophic to $S^n$, which in turn extends a result of Colding, Ilmanen, Minicozzi, and White \cite{CIMW} which says closed self shrinkers with entropy less than $\Lambda_{n-1} \lneq \Lambda_{n-2}$ are diffeomorphic to $S^n$, hence weakening the assumed entropy bound. In a similar manner, the result above extends (in a weaker sense than the compact case) a result of Bernstein and L. Wang \cite{BW1} for noncompact shrinkers in $\R^4$, where they showed (amongst other results, see corollary 1.4 therein) for asymptotically conical self shrinkers $M^3 \subset \R^4$ satisfying $\lambda(M) \leq \Lambda_2$ the stronger conclusion that they are diffeomorphic to $\R^3$; our argument does at least recover their statement though as discussed at the end of the proof. With the round cylinder in mind our conclusion seems likely to be sharp in this sense, although it could be possible that a shrinker in $\R^4$ with this entropy bound has more than one end precisely when it is a cylinder. 
 $\medskip$

In this dimension and under this entropy bound, we remark that generic mean curvature flow through neckpinch singularities has been established in \cite{CCMS, CCMS2}, so for some applications of the flow (see for instance \cite{DH}) the study of self shrinkers in this regime is unnecessary. However besides its intrinsic interest this result might still be of use in understanding singularity along nongeneric flows, which could imaginably occur for instance in problems involving families of flows (although to the author's knowledge potential fattening is a more serious concern). It also paints a picture explicitly how a perturbation of a nongeneric flow might only develop neckpinch singularities, by some copies of the $S^{2} \times \R_+$ in the statement above pinching off before roughly speaking the $S^3$ factor collapses to a point (as opposed to handles prematurely pinching off a more complicated model).
$\medskip$

An important extra difficulty to consider in the noncompact case versus the closed case is that apriori nontrivial topology may be ``lost'' to spatial infinity under the flow without being properly understood. To illustrate this concern by an admittedly crude thought experiment, a hypothetical translator asymptotically modeled on $T^2 \times \R$ would never develop a singularity and hence its topology would never be ``encountered'' as a high curvature region in the flow. In particular it seems for $n = 3$ asymptotically conical self shrinkers could apriori have a complicated link. Our first task, and really most of the work of this paper, will be to show that in fact the link is simple. 
\begin{thm}\label{thm2} Suppose $M^3 \subset \R^{4}$ is a smooth 2-sided asymptotically conical self shrinker with entropy less than $\Lambda_{1}$. Then its link $L$ is homeomorphic to a union of $S^{2}$. 
\end{thm} 
As an indication of why one might argue this is reasonable, consider that in general the link $L$ of a shrinker $M^n \subset \R^{n+1}$ is of dimension $n-1$ so an entropy bound of $\Lambda_{n-2}$ on $M^n$ implies morally that its link is low entropy: for a submanifold $N^k \subset \R^{k+1}$, we say $N$ is low entropy if $\lambda(N) < \Lambda_{k-1}$ (hence the title of paper, since $\Lambda_{k-2} > \Lambda_{k-1}$ as discussed in the next section). These compact surfaces with this entropy bound, at least in low dimensions ($n=2,3$) are known to be spheres. 
$\medskip$

The dimension bound assumption is for topological reasons (in the argument the classification of surfaces, Alexander's theorem, and Dehn's lemma are used which are dimension dependent, and in a probably less essential way the 3D Poincar\'e conjecture is also used) that perhaps indicate a deficit in knowledge and finesse moreso any true difficulty. As what one could take as a sliver of hope this is the case, Ilmanen and White in \cite{IW} showed lower bounds for the densities of area minimizing cones in terms of the topology of their link in every dimesnion. The area minimizing property there is employed by using a foliation near the cone by minimal surfaces, which are used as barriers in a mean curvature flow argument (at a high level our argument is similar to theirs). Since cones are noncompact, this is clearly the same sort of result as the topic of this paper.
$\medskip$

For instance, one might be hopeful to directly modify our argument in the next higher dimension ($n = 4$) because simply connected 3 manifolds are spherical by the resolution of the 3D Poincar\'e conjecture by Perelman \cite{P1, P2, P3} -- we use the corresponding (much easier) fact for surfaces below to classify the link. As an example of why this simple criterion alone doesn't seem to immediately lead to a proof of the corresponding statement for $n=4$, a potential issue (to the author's understanding) in this dimension is that the link could be a nontrivial homology sphere -- below we use that nonspherical oriented surfaces have nontrivial homology in a seemingly essential way. For higher dimensions of course there are higher dimensional versions of the Poincar\'e conjecture as verfied by Freedman and Smale \cite{Free, Smale}; this seems to naturally be even more complicated for a number of reasons than the $n=4$ case just discussed though. 
$\medskip$

$\textbf{Acknowledgements:}$ The author is supported by an AMS-Simons travel grant and thanks them for their generosity, as well as the anonymous referee for their careful reading and critique.

\section{Preliminaries}

Let $X: M \to N^{n+1}$ be an embedding of $M$ realizing it as a smooth closed hypersurface of $N$, which by abuse of notation we also refer to as $M$. Then the mean curvature flow $M_t$ of $M$ is given by (the image of) $X: M \times [0,T) \to N^{n+1}$ satisfying the following, where $\nu$ is the outward normal: 
\begin{equation}\label{MCF equation}
\frac{dX}{dt} = \vec{H} = -H \nu, \text{ } X(M, 0) = X(M)
\end{equation}
$\medskip$

By the comprison principle singularities occur often which makes their study important. To study these singularities, one may parabolically rescale about the developing high curvature region to obtain an ancient flow defined for times $(-\infty, T]$; when the base point is fixed this is called a $\textit{tangent flow blowup}$ which will be modeled on self shrinkers: by Huisken monotonicity \cite{H} these are surfaces satisfying the following equivalent definitions (at least when the shrinker is smooth but some definitions apply in the varifold sense as well): 
 \begin{enumerate}
\item $M^n \subset \R^{n+1}$ which satisfy $H - \frac{\langle X, \nu \rangle}{2} = 0$, where $X$ is the position vector 
\item  minimal surfaces in the Gaussian metric $G_{ij} = e^{\frac{-|x|^2}{2n}} \delta_{ij}$
\item surfaces $M$ which give rise to ancient flows $M_t$ that move by dilations by setting $M_t = \sqrt{-t} M$
\end{enumerate}

(These notions all make sense at least when the shrinker is smooth but some definitions apply in the varifold sense as well): As is well known, the second variation for formula for area shows there are no stable minimal surfaces in Ricci positive manifolds, see for instance chapter 1 of \cite{CM}. This turns out to also be true for minimal surfaces of polynomial volume growth in $\R^n$ endowed with the Gaussian metric as discussed in \cite{CM1}. To see why this is so, the Jacobi operator for the Gaussian metric is given by: 
\begin{equation}
L = \Delta + |A|^2 - \frac{1}{2} \langle X, \nabla(\cdot) \rangle + \frac{1}{2}
\end{equation} 
The extra $\frac{1}{2}$ term is essentially the reason such self shrinkers unstable in the Gaussian metric: for example owing to the constant term its clear in the compact case from this that one could simply plug in the function ``1'' to get a variation with $Lu >0$ which doesn't change sign implying the first eigenvalue is negative. In fact, every properly embedded shrinker has polynomial volume growth by Q. Ding and Y.L. Xin: 
\begin{thm}[Theorem 1.1 of \cite{DX}]\label{proper} Any complete non-compact properly immersed self-shrinker $M^n$ in $\R^{n+m}$ has Euclidean volume growth at most. \end{thm}
We will combine these facts below to conclude the self shrinker we could find in some cases in the argument below must in fact be unstable. 
$\medskip$

The mean curvature flow is best understood in the mean convex case because it turns out under quite weak assumptions the only possible shrinkers are generalized cylinders $S^k \times \R^{n-k}$, especially so for 2-convex surfaces ($\lambda_1 + \lambda_2 > 0$) and a surgery theory with this convexity condition similar to the Ricci flow with surgery has been carried out. For the mean curvature flow with surgery one finds for a 2-convex surface $M$ curvature scales $H_{th} < H_{neck} < H_{trig}$ so that when $H = H_{trig}$ at some point $p$ and time $t$, the flow is stopped and suitable points where $H \sim H_{neck}$ are found to do surgery where ``necks'' (at these points the surface will be approximately cylindrical) are cut and caps are glued in. The high curvature regions are topologically identified as $S^n$ or $S^{n-1} \times S^1$ and discarded and the low curvature regions will have curvature bounded on the order of $H_{th}$. The flow is then restarted and the process repeated. 
$\medskip$

It was initially established for compact 2-convex hypersurfaces in $\R^{n+1}$ where $n \geq 3$ by Huisken and Sinestrari in \cite{HS}, and their approach was later extended to the case $n =2$ by Brendle and Huisken in \cite{BH}, where 2-convexity is mean convexity. A somewhat different approach covering all dimensions simultaneously was given later by Haslhofer and Kleiner in \cite{HK} shortly afterwards. Haslhofer and Ketover then showed several years later in section 8 of their paper \cite{HKet} enroute to proving their main result that the mean curvature flow with surgery can be applied to \textit{compact} mean convex hypersurfaces in general ambient manifolds. Important to this article, the author with S. Wang established it for (compact) mean convex hypersurfaces with entropy less than $\Lambda_{n-2}$ in the sense of Colding and Minicozzi:
$\medskip$

In \cite{CM}, Colding and Minicozzi introduced their important notion of entropy, which is defined as the supremum of translated and rescaled Gaussian densities; indeed, consider a hypersurface $\Sigma^k \subset \R^{\ell}$; then given $x_0 \in \R^{\ell}$ and $r > 0$ define the functional $F_{x_0, r}$ by 
\begin{equation}
F_{x_0, r}(\Sigma) = \frac{1}{(4 \pi r)^{k/2}} \int_\Sigma e^\frac{-|x - x_0|^2}{4r} d\mu
\end{equation} 
(when $x_0 = \vec{0}$ and $r = 1$, this is just a normalization of area in the Gaussian metric). Colding and Minicozzi then define the entropy $\lambda(\Sigma)$ of a submanifold to be the supremum over all $F_{x_0, r}$ functionals:
\begin{equation} 
\lambda(\Sigma) = \sup\limits_{x_0, r} F_{x_0, r}(\Sigma) 
\end{equation} 
The aforementioned Huisken monotonicity \cite{H} implies that this quantity is in fact monotone under the flow, and because it is defined as a supremum over rescalings and recenterings it also controls the nature of singularities encountered along the flow -- see \cite{CIMW, BW0, BW, BW1} for instance. Note that surfaces of polynomial volume growth have finite entropy. 
$\medskip$

The current state of knowledge of mean curvature flow singularities approached from an entropy perspective seems to be ``quantized'' by the entropy $\Lambda_k$ of round spheres as we now discuss. By a calculation of Stone \cite{St} we have: $$\Lambda_1>\frac{3}{2}>\Lambda_2>...>\Lambda_n\rightarrow\sqrt2$$
So far in the literature, many results using an entropy condition assume that the submanifold $M$ under consideration satisfies $\lambda(M) < \Lambda_{n-1}$, which seem to most often be refered to as a low or small entropy condition. The next natural entropy condition to consider then is a bound by $\Lambda_{n-2}$, which in this paper we refer to naturally as a medium entropy bound; one might expect studying surfaces with this entropy bound to be tractable because morally it implies that mean convex singularities encountered will be $2$-convex, which as implied above in the discussion on surgery are the easiest to consider/flow through (after convex ones). Indeed this philosophy was carried out in the compact case in the the joint work with S. Wang \cite{MW} (``low'' in its title refers to what we define as medium). An important observation for the argument of this paper is that this philosophy can be extended to the noncompact setting, but there are some important new issues to consider. For instance in the noncompact case the asymptotics of the submanifold in question are important: 
$\medskip$

Throughout this article we will say an end $E$ of a self shrinker is \textit{asymptotically conical} if $E$ satisfies $\lim_{\tau \to \infty} \tau^{-1} E = C(E)$ in $C_{loc}^\infty(\R^{n+1} \setminus 0)$ for $C(E)$ a regular cone in $\R^{n+1}$. A similar definition can be made for asymptotically cylindrical ends and by results of L. Wang \cite{Lu} for $n = 2$ every end of a self shrinker of finite topology is either asymptotically conical or cylindrical (with multiplicity one). Naturally, one says a self shrinker is asymptotically conical if every end is. Considering singular/GMT extensions of shrinker and asymptoticaly conical end in a natural way, under suitable entropy assumptions and assumptions of the underlying measure the support of the shrinker and asymptotic cone can be shown to be smooth (so asymptotically conical as in the sense above) -- see propositions 3.2, 3.3 in \cite{BW1} and lemma 2.1 in \cite{CCMS2}. In particular, our theorem applies to asymptotically conical (in the weak sense) shrinkers which arise as blowups under our entropy assumption in $\R^4$. Note, since $\Lambda_1 < 2$, the convergence will be with multiplicity one.   
$\medskip$

In the same paper where they introduced entropy, Colding and Minicozzi showed the only singularities which morally shouldn't be able to perturbed away are the mean convex ones, the generalized round cylinders $S^k(\sqrt{2k}) \times \R^{n-k}$, called round because the spherical factor is a standard round sphere of a radius appropriate to satisfy the shrinker equation. In particular other singularity models should be able to be perturbed away, so round cylinders are called generic singularity models. Their numbers are few (only $n$ of them), wheras for instance in $\R^3$ there are self shrinkers are of arbitrarily large genus by \cite{KKM}, so one could say most self shrinkers are nongeneric. 
$\medskip$

Concerning nongeneric singularity models, the no-cylinder conjecture of Ilmanen says that the types of ends shouldn't be ``mixed'' in that if there is a single cylindrical end then $M$ is a cylinder so one expects that ``most'' self shrinkers in $\R^3$ are asymptotically conical (see \cite{Lu1} for a partial result confirming this). Extending this conjecture to the next higher dimension, this provides our justification with the above paragraph in mind for the claim that self shrinkers are ``often'' asymptotically conical -- it is also quite convenient for analytical reasons. 
$\medskip$

Returning to flows through singularities, an important advantage of the mean curvature flow with surgery is that the topological change across discontinous times, when necks are cut and high curvature regions discarded, is easy to understand. A disadvantage is that it isn't quite a Brakke flow (a geometric measure theory formulation of the mean curvature flow) so does not immediately inherit some of the consequences thereof, but at least it is closely related to the level set flow by results of Laurer \cite{Lau} and Head \cite{Head, Head1} which in the nonfattening case is (modulo some technicalities). In their work they show that surgery converges to the level set flow in Hausdorff distance (and in fact in varifold sense as Head shows) as the surgery parameters degenerate (i.e. as one lets $H_{th} \to \infty$). This connection is useful for us because deep results of White \cite{W} show that a mean convex LSF will converge to a (possibly empty) stable minimal surface long term. 
$\medskip$

As mentioned above mean curvature flow with surgery in a curved ambient setting (at least for 3 manifolds and bounded geometry) has been already accomplished by Haslhofer and Ketover but some extra care is needed for the Gaussian metric especially in the noncompact case. This is because the metric is poorly behaved at infinity (as one sees from the calculation of its scalar curvature) which introduces some analytic difficulites for using the flow, so instead we consider the renormalized mean curvature flow (which we'll abreviate RMCF) defined by the following equation: 

\begin{equation}\label{renorm1}
\frac{dX}{dt} =  \vec{H} + \frac{X}{2}
\end{equation}

Where here as before $X$ is the position vector on $M$. It is related to the regular mean curvature flow by the following reparameterization; this allows one to transfer many deep theorems on the MCF to the RMCF. Supposing that $M_t$ is a mean curvature flow on $[-1,T)$, $-1 < T \leq 0$ ($T = 0$ is the case for a self shrinker). Then the renormalized flow $\hat{M_\tau}$ of $M_t$ defined on $[0, -\log(-T))$ is given by 
\begin{equation}\label{param}
\hat{X}_\tau = e^{\tau/2} X_{-e^{-\tau}},\text{ } \tau = -\log{(-t)}
\end{equation}

Note up to any finite time the reparameterization is bounded and preserves many properties of the regular MCF, like the avoidance principle and that entropy is monotone under the RMCF. With this in mind, the author showed in his previous article \cite{Mra1} that one can then construct a flow with surgery using the RMCF on suitable perturbations of noncompact self shrinkers, and that as one lets the surgery parameters degenerate indeed the surgery converges to the level set flow when $n =2$. This can be readily combined with the aforementioned joint work with S. Wang \cite{MW} to show the following:

\begin{thm}\label{LSF} Let $M^n \subset \R^{n+1}$ be a smoothly asymptotically conical hypersurface such that $H - \frac{\langle X, \nu \rangle}{2} \geq c(1 + |X|^2)^{-\alpha}$ for some constants $c, \alpha >0$ and choice of normal, that $\lambda(M) < \Lambda_{n-2}$. Then denoting by $K$ the region bounded by $M$ whose outward normal corresponds to the choice of normal on $M$, the level set flow $M_t$ of $M$ with respect to the renormalized mean curvature flow satisfies
\begin{enumerate}
\item inward in that $K_{t_1} \subset K_{t_2}$ for any $t_1 > t_2$, considering the corresponding motion of $K$.
\item $M_t$ is the Hausdorff limit of surgery flows $S_t^k$ with initial data $M$. 
\item $M_t$ is a forced Brakke flow (with forcing term given by position vector). 
\end{enumerate}
\end{thm}

$\alpha$-noncollapsedness here means there are inner and outer osculating balls of radius proportional to the shrinker mean curvature and this has many consequences, see \cite{BA,HK0}. The assumption on the asymptotics are conditions for which shrinker mean convexity are preserved and existence of a entropy decreasing perturbation of a self shrinker smoothly asymptotic to a cone can always be assumed to satisfy this by work of Bernstein and L. Wang in \cite{BW}. We use this theorem (often implicitly) below with $n = 3$ when we discuss the flow of $M$. 
$\medskip$

$H_G$, the mean curvature in the Gaussian metric, is related to the renormalized mean curvature by $H_G = e^{\frac{|x|^2}{4} }(H - \frac{\langle X, \nu \rangle}{2})$ and as a result the time limit of the flow defined in the theorem above by White's theory for mean convex MCF (in particular \cite{W}) will be a stable self shrinker if nonempty. It will also have finite entropy by Huisken monotonicity and hence have polynomial volume growth. As a result either by the instability results mentioned above or by the Frenkel theorem for self shrinkers given in the appendix of \cite{CCMS} we have the following:
\begin{lem}\label{Wthy} Let $M_t$ be the flow defined above in theorem \ref{LSF}. Then $\lim\limits_{t \to \infty} M_t = \emptyset$. 
\end{lem}

Lastly, note that by switching our choice of normal and using minimality (in Gaussian metric) of the original surface we may shrinker mean convex perturb either inward or outward (for a 2-sided surface, of course the distinction is somewhat arbitrary) to study its topology as observed in \cite{HW, BW} -- this idea is critical to our argument and we will choose our choice of perturbation depending on which case we are considering in the argument below.

\section{Proof of Theorem \ref{thm2}} 

Note by the entropy assumption (in particular that $\Lambda_1 < 2$) that $M$ is embedded and hence its link is too. Without loss of generality for this section, the link $L$ is connected. Supposing $L$ is not diffeomorphic to $S^{2}$, there exists some $R_1 >> 0$ so that $M \cap S(0,R) := L_R$ for is not diffeomorphic to $S^{2}$ for $R > R_1$ and that, by the asymptotically conical assumption, $L_R \simeq L_{R_1}$ for all $R > R_1$. By the classification of surfaces (note $L_R$ is orientable, which can be seen by projecting the normal of $M \cap S(0,R)$ onto $TS(0,R)$, giving a section of the normal bundle of $L_R$, which has no kernel because the sphere intersects $M$ transversely) $L_R$ is topologically a connect sum of tori which bounds domains (not necessarily handlebodies) $K_R, K_R^c \subset S(0,R)$. 
$\medskip$

 Fixing a choice of $R > R_1$, consider a standard generator $\gamma \subset L_R$ of $H_1(L_R)$ (as in, writing $L_R$ as a connect sum of tori, $\gamma$ is homotopic to one of the two generators of a single one of the tori). Note $\gamma$ is also homotopically nontrivial in $L_R$. We consider two cases: either $\gamma$ is homotopically trivial in $M$ or not. Without loss of generality, $\gamma$ is embedded and smooth as well. 
 $\medskip$
 \begin{figure}
\centering
\includegraphics[scale = .5]{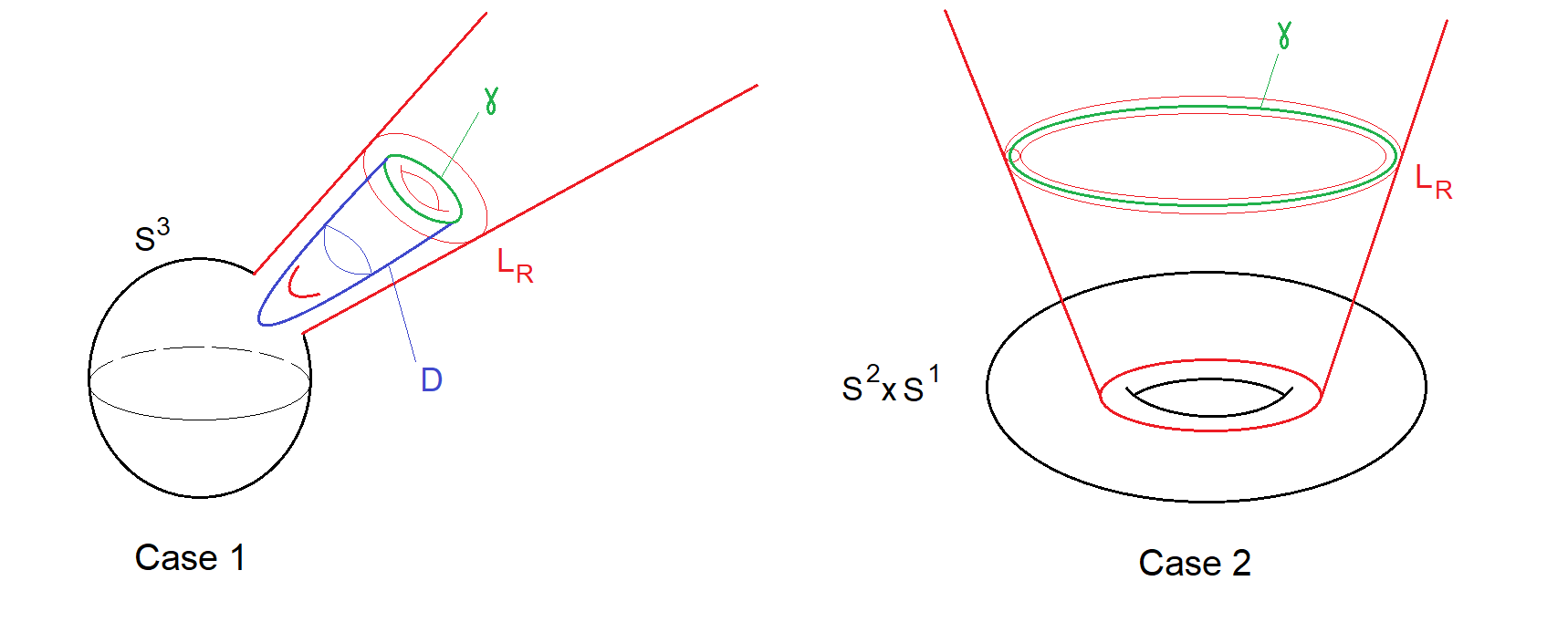}
\caption{This figure illustrates toy examples one might imagine for the two cases considered when $L_R$ is a standardly embedded torus. In the first case $\gamma \in \partial (M \cap B(0,R))$ is nullhomotopic in $M$ and hence bounds an embedded disc $D \subset M \cap B(0,R)$ by Dehn's lemma. }
\end{figure}

\subsection{First case: $\gamma$ is nullhomotopic in $M$}
$\medskip$

Since $\gamma$ is homotopically trivial in $M$, it bounds a disc $D$ in $M$; suppose $D \subset M \cap B(0, R_2)$. Hence for any embedded curve $\gamma' \subset L_{R'}$ isotopic (in $M$) to $\gamma$ for $R' > R_2$, $\gamma'$ is nullhomotopic in $M \cap B(0, R')$ and hence bounds an embedded disc $D' \subset M \cap B(0,R')$ by Dehn's lemma (see \cite{Hat})  -- Dehn's lemma gives a PL embedded disc but when $\gamma'$ is smooth note $D'$ can be taken to be smooth as well by the Whitney approximation theorem \cite{Lee}. The idea is morally such discs serve as barriers in a sense to keep the flow of (a perturbation of) $M$ ``propped'' up. The following indicates which domain $M$ bounds to perturb and flow into: 
$\medskip$

\begin{lem}\label{nontrivial} In one of $K_R$ or $K_R^c$ the curve $\gamma$ is not homotopically trivial. 
\end{lem}
\begin{pf}
It seems one could probably use the Mayer-Vietoris sequence and Hurwicz isomorophism here as in the proof of theorem 1 in \cite{HW} but we present a more geometric argument. Suppose for the sake of contradiction it were homotopically trivial in both simultaneously. By Dehn's lemma $\gamma$ bounds PL embedded (of course, in fact smooth) discs $D_1 \subset K_R$ and $D_2 \subset K_R^c$ which intersect along $\gamma$ giving an embedded $S^2 \subset S^3$. Since $\gamma$ is smooth their union gives a PL embedded $S^2$ and so by Alexander's theorem (see \cite{Hat}) $D_1 \cup D_2$ then bounds a (PL embedded) 3-ball $B \subset S^3$. From its construction $L_R$ intersects $B$ in one boundary component, namely $\gamma$. In particular, $\gamma$ is homologically trivial in $L_R$ giving a contradiction. 
\end{pf}
$\medskip$

Of course, the lemma applies equally for $\gamma'$ homotopic to $\gamma$ in $L_{R'}$ for $R' > R > R_1$. After potentially relabeling, $\gamma$ is homotopically nontrivial in $K_R$. In this case, consider a shrinker mean convex perturbation of $M$, as constructed in \cite{BW}, which descends (i.e. intersecting with $S(0,R)$) to a perturbation of $L_R$ into $K_R$ and consider the corresponding renormalized flow $M_t$ (recalling we can choose which direction to flow into as discussed in the preliminaries). This flow likewise descends to a flow $(L_R)_t$ of $L_R$. Note though that although $M_t$ is a RMCF that $(L_R)_t$ isn't necessarily (to the author's knowledge) an easily described flow in $S(0,R)$, but we will still find it profitable to consider. 
$\medskip$

 By lemma \ref{Wthy} $M_t$ must leave every bounded set in some finite time and hence $(L_R)_t$ must eventually become empty. Denote this time by $T$. We will play it off against the next two lemmas, the first essentially that the disc we find by Dehn's lemma persists: 

\begin{lem}\label{bound} For any $\tau > 0$, one can pick $\overline{R}$ sufficiently large so that for $t \in [0, \tau]$, there will be a smoothly embedded curve $\overline{\gamma}_t \subset (L_{\overline{R}})_t$ isotopic to $\gamma$ which will bound a smoothly embedded disc $\overline{D}_t \subset B(0,\overline{R}) \cap M_t$. 
\end{lem}

\begin{pf} 
 With the construction of the flow by surgery flows given in thoerem \ref{LSF} in mind, we first show for exposition that this holds for an approximating surgery flow $S_t$ of $M$. Since the $L_{R'}$ are all diffeomorophic for $R'$ large there is clearly an initial choice of curve $\overline{\gamma}$ which is isotopic to $\gamma$. Up to the first surgery time and in between surgery times when the surgery flow is smooth, this curve is just given by restricting the motion of $(L_{\overline{R}})_t$ along $\overline{\gamma}$ because for $\overline{R}$ large enough $(L_{\overline{R}})_t$ will be a graph over $(L_{\overline{R}})_0$ by pseudolocality and item (3) of theorem \ref{LSF} on $[0, \tau]$ (in particular, $\overline{\gamma}_t$ can be taken to be embedded and smoothly vary) for times in $[0, \tau]$. Concerning the bounded disc for smooth times the flow is an isotopy which restricts to an isotopy of the disc (modding out tangential components of the flow). Now consider a surgery time $t_s < \tau$ -- we must check that after surgery $\overline{\gamma}_{t_s}$ still bounds a disc. Again by pseudolocality for all surgery necks $N$, $N \cap S(0,\overline{R})$ is empty and similarly all $N$ must be within $B(0,\overline{R})$. Considering a cap $C$ in the surgery procedure, since it is topologically a ball the intersection of $D_{t_s}$ with $\partial C \simeq S^2$ is a disjoint union of closed curves which bound discs by the Schoenflies theorem (without loss of generality $D_t$ enters all caps transversely). Surgering along these discs gives a union of $S^2$ along with a new disc whose boundary is $\overline{\gamma}_{t_s}$ (essentially filling in the part of the disc between the end and the ``closest'' surgery necks).  In particular, $\overline{\gamma}_{t_s}$ continues to bound a disc after surgery. Note that its concievable at this stage that the discarded copies of $S^2$ bound nontrivial topology of $M_t$, so the $D_t$ do not necessarily form an isotopic family of discs apriori.
$\medskip$

Now we discuss how to show the curve $\overline{\gamma}$ from the previous paragraph always bounds a disc in the limiting flow. What one might first wish for is to take a limit (by compactness) of the discs as the surgery parameters degenerate but if the limiting disc enters a singular region of the flow it could potentially complicate things so its best if the disc is taken to avoid it completely. There is also the matter of boundedness along this sequence of discs needed to apply a compactness theorem, which suggests its best in terms of the disc to work only within the context of the level set flow. 
$\medskip$

To begin, we consider high curvature points we might encounter as we travel sufficiently deep within a high curvature region (loosely speaking) from a low curvature region as in our situation of a disc starting from an end (where $\overline{\gamma}$ is) approaching a singularity in the interior of $B(0,\overline{R})$. At points where say where $H \sim H_{can}$, refering to parameters in the canonical neighborhood theorem (see \cite{HK}, here we are supressing some notation as well), one can find nearby ``neck like'' points (see proposition 3.2 in \cite{HK}) in any corresponding ancient model that could appear irrespective of surgery parameters. Intuitively a surgery flow near such a point is modeled locally by a neck or a cap bordered by a neck facing towards the low curvature region. With the Hausdorff convergence in mind then, one can use Arzela-Ascoli to pass to the limit on these bounded curvature regions for the surgery flows to see that the level set flow always has necks where $H \sim H_{can}$ as one approaches a singular region of the level set flow from a low curvature one (of course, these are smooth points as well). These necks on the level set flow give points to surger the disc as we did for the surgery flows; in this case we perform the surgery on the disc whenever a point of it is in a region of the level set flow where say $H = 10 H_{can}$. Note between these times the disc varies continuously since it is within a region of the level set flow of bounded curvature, and that the disc can be taken to be smooth at all times since it is surgered on at along cross sections of necks of bounded curvature. 
\end{pf}

Without using Dehn's lemma (and in particular with replaying the argument in higher dimensions in mind), it seems the intersection of the disc with the boundary of the cap could be much more complicated although naively it seems likely that $\gamma$ remains homotopically trivial will still hold. We will pit this lemma against the defintion of the time $T$ by the following lemma, which says the discs must  ``leave'' the end no matter what: 
$\medskip$

\begin{lem}\label{leave} With $\overline{\gamma}$ and $\tau$ as in lemma \ref{bound}, after potentially taking $R$ larger, there is an $\epsilon > 0$ so that in $B(0, R - \epsilon)^c \cap M_t$ the curve $\overline{\gamma}_t$ isn't nullhomotopic and so doesn't bound a disc $B(0, R - \epsilon)^c \cap M_t$. In particular, the disc $D_t$ from the previous lemma satisfies that $D_t \cap S(0,R - \epsilon)$ is nonempty on $[0, \tau]$.
\end{lem}

\begin{pf}
 Denote by $K$ to be the region $M$ bounds which includes $K_R$. Note then that in $B(0, R - \epsilon)^c \cap K \simeq K_R \times [R -\epsilon, \infty)$ (for $R$ large enough), and in particular $\overline{\gamma}$ is homotopically nontrivial in this domain since it is homotopically nontrivial in $K_R$. If for some time $t \in [0, \tau]$ $ \overline{\gamma}_t$ is nullhomotopic in $B(0, R - \epsilon)^c \cap M_t$, then in particular $\overline{\gamma}_t$ bounds (the image of) a disc in $B(0, R - \epsilon)^c \cap K_t$. By the set monotonicity of the flow i.e. that $K_t \subset K$ we get in fact $\overline{\gamma}_t$ and hence $\overline{\gamma}$ are nullhomotopic in $B(0, R - \epsilon)^c \cap K$, giving a contradiction. 
\end{pf} 
$\medskip$

Applying the above lemmas with $\tau = T + 1$ we see we arrive at a contradiction. Considering a time $t \in (T, T+ 1)$ and the disc $D_t$ given from lemma \ref{bound}, the disc by lemma \ref{leave} must have nonempty intersection with $S(0, R - \epsilon)$. On the other hand, it cannot pass through $S(0,R)$ because $(L_R)_t = \emptyset$ for $t > T$. This completes the argument in this case.

\subsection{Second case: $\gamma$ homotopically nontrivial in $M$}
$\medskip$

Now consider the second case, that $\gamma$ is homotopically nontrivial in $M$ -- this case is easier in a sense because we may directly apply the deep ideas of White \cite{W3}, in particular theorems 1.1 and 5.2 therein. Specialized to our setting, a consequence of it is that if $K$ is a smooth mean convex set (and compact as stated in theorem 1.1, but this can also apply in the noncompact case as long as singularities occur only in a bounded ball for the time it is applied by theorem 5.2) in a Riemannian manifold $N$ of dimension 4 (in particular, less than 7), then if a curve in $K^c$ is initially homotopically nontrivial and later becomes contractible in $(K^c)_t$ a singularity of the form $S^1 \times S^2$ must have occured, contradicting the entropy bound. Here we will consider $N$ to be a subset of $\R^4$ (possibly all of $\R^4$ depending on which case we are in below) endowed with the Gaussian metric so the flow constructed is more precisely a mean convex foliation; however the flow is monotone, satisfies the Brakke regularity theorem and the singular set dimension results of White \cite{W0}, and all the singularities are modeled on round cylinders so the results of the paper apply -- this is essentially the upshot of Hershkovits and White \cite{HW} (although they phrase things enitrely in terms of the RMCF) where they study the interplay of entropy and topology for compact self shrinkers. 
$\medskip$

There are two possible cases for $\gamma$: that $\gamma$ is homotopically nontrivial in one of the components $K$, $K^c$ of $\R^4$ it bounds or not. First suppose suppose that $\gamma$ is homotopically nontrivial in (at least) one of $K$ or $K^c$, say $K^c$ to allign with White's terminology. Consider then a nontrivial curve $\gamma$ in $K^c$. Since $\gamma$ is contractible in $\R^4$, the corresponding homotopy gives it bounds a (continuous image of, perhaps not embedded) disc $D$ -- note this disc must intersect $K$. Perturbing and flowing into $K$ by lemma \ref{Wthy} eventually we must have $D \subset K^c_t$, say by $T$, implying by this time that $\gamma'$ is nullhomotopic in $K^c_t$. By psuedolocality \cite{CY} there is an $R >> 0$ so that near $S(0,R)$ $M_t$ is a smooth flow which intersects the sphere transversely, so defining $N = B(0,R)$ section 5 of White \cite{W3} implies a singularity modeled on $S^1 \times \R^2$ formed contradicting the entropy bound. 
$\medskip$

Now we consider the possibility that $\gamma$ is homotopically trivial in both $K$ and $K^c$; this naively seems to be a more exotic case than above, but we are unsure it can be ruled out apriori by purely topological reasoning. Then $\gamma \in M$ bounds a disc in both $K$ and $K^c$. Picking essentially arbitrarily (only to align with White's notation), we define $\widetilde{N}$ to be the union of $K$ and $K^c \cap M \times [0, \epsilon)\nu$ (i.e. a collar of $M$), where $\nu$ is the normal pointing away from $K$ and $\epsilon > 0$ is some number small enough so that the $\epsilon$ level set of the collar is also embedded in $\R^4$. Note this collar region retracts onto $M$; the utility of this is that now $\gamma$ is a homotopically nontrivial curve in $K^c \cap \widetilde{N} \subset \widetilde{N} \subsetneq \R^4$. Consider as in the previous paragraph a disc $D \subset K$ bounding $\gamma$ and flow out of $K^c$ into $K$ (that the disc can be taken to be contained in a single component, and hence in $N$, is why we split this up into cases). As above, lemma \ref{Wthy} (this still applies since the flow of $M_t$ is the same considered in $\R^4$ or $\widetilde{N}$) gives that eventually $D \subset K^c_t \cap \widetilde{N}$ -- call this time $T$. Let $R >> 0$ be large enough so that $M_t$ intersects $S(0,R)$ only transversely and as a smooth flow; again such an $R$ exists by pseudolocality. Defining $N = \widetilde{N} \cap B(0,R)$ then and noting $\gamma$ is still homologically nontrivial in $K^c \cap \widetilde{N} \cap B(0,R)$, section 5 of \cite{W3} give that a singularity modeled on $S^1 \times \R^2$ must have formed (in fact by time $T$) giving a contradiction. 
$\medskip$

\section{Proof of theorem \ref{thm1}}

By the Frenkel property for self shrinkers, $M$ must be connected. By theorem \ref{thm2} then, there exists $R$ sufficiently large so that $M \cap B(0,R)$ is diffeomorphic to a connected, 2-sided hypersurface $N^3$ whose boundary consists of a number of 2-spheres along each of which an end homeomorphic to $S^2 \times \R_+$ is attached, where by ends here we mean for an appropriate choice of $R$ disjoint connected components of $M \setminus B(0, R)$ which are diffeomrophic to half cylinders over distinct (for distinct ends) connected components of the link; such an $R$ exists since $M$ is asymptotically conical, and the convergence is multiplicity one. The point is to confirm that $N$ is simply connected; then by capping off each component of $N$ (considering $N$ as an intrinsically defined manifold, as a hypersurface in $\R^4$ it seems some ends could be ``parallel'' which would preclude doing this at least in an embedded way) with a 3-ball we obtain a closed, connected, simply connected 3 manifold $\widetilde{N}$ which by the resolution of the 3D Poincar\'e conjecture is diffeomorphic to $S^3$. If there is a homotopically nontrivial curve on $N$ and hence $M$ by the Seifert-Van Kampen theorem and theorem \ref{thm2}, we can proceed directly as we do in the second case of the proof above using \cite{W3}, giving the first part of theorem \ref{thm1}. Note that with surgery for compact manifolds in mind one should be able to argue directly with a bit more work that $\widetilde{N}$ is diffeomoroprhic to either $S^3$ or a connect sum of $S^2 \times S^1$, the later of which could subsequently be ruled out avoiding the use of the Poincar\'e conjecture -- this seems to be naturally a more robust line of reasoning for considering higher dimensional versions of our statement. 
$\medskip$

When the number of ends is equal to one, $M$ is diffeomorphic to $\R^3$ as a consequence of Alexander's theorem as noted in \cite{BW1}. Now suppose that $\lambda(M) < \Lambda_2$ (we will discuss the case of equality afterwards) and $M$ had (at least) two ends, labeled $E_1$ and $E_2$. Fixing an $R$ in our definition of end given in the paragraph above, consider a curve $\gamma: \R \to M$ such that for $s$ sufficiently negative $\gamma(s) \in E_1$ and $s$ sufficiently positive $\gamma(s) \in E_2$. With this in mind intersect $M$ with an embedded hypersurface $P \simeq \R^3$ such that (i) $E_1$ lays on one side of $P$ and $E_2$ lays on the other side, (ii) $P$ intersects $M$ transversely, and (iii) that $P \cap M$ is compact; this is always possible by the asymptotically conical assumption. Denote by $P \cap M$ the surface $S$; note that $S$ is closed since its compact and $M$ is boundaryless. Similarly denote the bounded portion of $P$ that $S$ bounds by $K_S$. By perturbing and flowing $M$ so that, restricted to $N$, the flow is into $K_S$ (using (ii)) we see as above by lemma \ref{Wthy} that $S_t$ is eventually empty. Because of this one may argue that a singularity of $M_t$ must occur which disconnects $E_1$ from $E_2$ along $\gamma$; note that by using large spheres as barriers far along the ends toward spatial infinity (or, alternatively, pseudolocality), for any given finite time there will be points orginating from $E_1$ and $E_2$ on one side of $P$ and the other, respectively, or in other words one end can't flow from one side of $P$ to the other in finite time so a singularity which disconnects $M$ must indeed occur. Clearly such a singularity must be modeled on $S^2 \times \R$, which has entropy $\Lambda_2$ contradicting $\lambda(M) < \Lambda_2$. In the case $\lambda(M) = \Lambda_2$, we note the perturbation of Bernstein and Wang we used strictly decreases entropy placing us into the case of strict inequality.

\end{document}